\documentclass[12pt,a4paper]{article}
\usepackage{bbm}
\usepackage{amsmath,amssymb}
\usepackage{setspace}
\usepackage{algpseudocode,algorithm}
\usepackage{tabularx}
\usepackage{amsthm}
\usepackage{graphicx}
\usepackage{here}
\usepackage{url}
\usepackage{txfonts}

\usepackage{natbib}
\bibpunct[:]{(}{)}{;}{a}{}{;}
\usepackage[usenames]{color}
\usepackage{authblk}

\allowdisplaybreaks[1]

\newtheorem{theorem}{Theorem}[section]
\newcommand{\argmin}{\mathop{\rm argmin}\limits}
\makeatletter
  
  \@addtoreset{equation}{section}
\makeatother

\def\gd{$\gamma$-divergence}
\def\pdf{probability density function}
\def\eps{\varepsilon}
\def\iidp{i.i.d.\,problem}
\usepackage{bm}

\title{On Difference Between Two Types of $\gamma$-divergence for Regression}
\date{} 

\author[1]{Takayuki Kawashima} 
\author[1,2,3]{Hironori Fujisawa} 

\affil[1]{ Department of Statistical Science, The Graduate University for Advanced Studies, Tokyo  \authorcr  E-mail: \texttt{t-kawa@ism.ac.jp}}

\affil[2]{The Institute of Statistical Mathematics, Tokyo  \authorcr E-mail: \texttt{fujisawa@ism.ac.jp} }

\affil[3]{ Department of Mathematical Statistics, Nagoya University Graduate School of Medicine  }

\begin{document}
\maketitle

\begin{abstract}
The $\gamma$-divergence is well-known for having strong robustness against heavy contamination. By virtue of this property, many applications via the $\gamma$-divergence have been proposed. There are two types of \gd\ for regression problem, in which the treatments of base measure are different. In this paper, we compare them and pointed out a distinct difference between these two divergences under heterogeneous contamination where the outlier ratio depends on the explanatory variable. One divergence has the strong robustness under heterogeneous contamination. The other does not have in general, but has when the parametric model of the response variable belongs to a location-scale family in which the scale does not depend on the explanatory variables or under homogeneous contamination where the outlier ratio does not depend on the explanatory variable. \citet{hung.etal.2017} discussed the strong robustness in a logistic regression model with an additional assumption that the tuning parameter $\gamma$ is sufficiently large. The results obtained in this paper hold for any parametric model without such an additional assumption. 
\end{abstract}

\section{Introduction}\label{Sect: Intro}

The maximum likelihood estimation and least square method have been widely used. However, they are not robust against outliers. To overcome this problem, many robust methods have been proposed mainly using M-estimation \citep{hampel.etal.2005,maronna:martin:yohai:2006,huber.ronchetti.2009}. The maximum likelihood estimation can be regarded as the minimization of the empirical estimator of the Kullback-Leibler divergence. As an extension of this idea, some robust estimators were proposed as the minimization of the empirical estimators of the modified divergences \citep{basu.etal.1998,scott.2001,fujisawa.eguchi.2008}. 

Recently, some robust regression methods based on divergences have been proposed using $L_2$-divergence \citep{chiandscott2014,lozano2016}, density power divergence \citep{Ghosh2016} and $\gamma$-divergence \citep{hung.etal.2017,kawashim.fujisawa.2017}. The robust properties are often investigated under the contaminated model. The difference between i.i.d. and regression problems is whether the outlier ratio in the contaminated model can depend on the explanatory variable or not. They are called the heterogeneous and homogeneous contaminations, respectively. To our knowledge, the robust properties under heterogeneous contamination have not been investigated well. Recently, \citet{hung.etal.2017} pointed out that a logistic regression model with mislabel can be regarded as a logistic regression model with heterogeneous contamination and then applied the \gd\ to a usual logistic regression model, which enables us to estimate the parameter of the logistic regression model without estimating a mislabel model even if mislabels exist. They discussed the strong robustness that the latent bias can be sufficiently small against heavy contamination, but assumed that the tuning parameter $\gamma$ is sufficiently large.

There are two types of \gd\ for regression problem in which the treatments of base measure are different \citep{fujisawa.eguchi.2008,kawashim.fujisawa.2017}. In this paper, we compare them in detail and prove that one divergence can show the strong robustness for any parametric model under heterogeneous contamination. The other can not in general, but can when the parametric model of the response variable belongs to a location-scale family in which the scale does not depend on the explanatory variables or under homogeneous contamination. This difference will be illustrated in numerical experiments. Comparing them with Hung et al. (2017), the results obtained here holds for any parametric model, including a logistic regression model, and do not assume that $\gamma$ is sufficiently large. 

This paper is organized as follows. 
In Section~\ref{Sect: reg based gam-div}, we review two types of $\gamma$-divergence for regression problem.
In Section~\ref{Sect: rob prop}, we elucidate a distinct difference between two types of $\gamma$-divergences from the viewpoint of robustness. 
In Section~\ref{Sect: experment}, numerical experiments are illustrated to verify the difference stated in Section~\ref{Sect: rob prop}. 
The R language code, which was run in our experiments, is available at \url{https://sites.google.com/site/takayukikawashimaspage/software}.

\section{Regression based on $\gamma$-divergence}\label{Sect: reg based gam-div}
The $\gamma$-divergence for regression was first proposed by \citet{fujisawa.eguchi.2008}. 
It measures the difference between two conditional probability density functions.
The other type of the $\gamma$-divergence for regression was proposed by \citet{kawashim.fujisawa.2017}, in which the treatment of the base measure on the explanatory variable was changed. In this section, we briefly review both types of $\gamma$-divergences for regression and present the corresponding parameter estimation.

\subsection{Two types of $\gamma$-divergences for regression}\label{Sect def two gamma}

First we review the \gd\ for the \iidp. Let $g(u)$ and $f(u)$ be two probability density functions. The $\gamma$-cross entropy and \gd\ were defined by
\begin{eqnarray*}
d_\gamma(g(u),f(u)) &=& -\frac{1}{\gamma} \log \int g(u) f(u)^\gamma du + \frac{1}{1+\gamma} \log \int f(u)^{1+\gamma} du, \\
D_\gamma(g(u),f(u)) &=& - d_\gamma(g(u),g(u)) + d_\gamma(g(u),f(u)) .
\end{eqnarray*}
This satisfies the following two basic properties of divergence: 
\begin{center}
\begin{tabular}{cl}
(i) & $D_\gamma(g(u),f(u)) \ge 0$. \\
(ii) & $D_\gamma(g(u),f(u)) = 0 \ \Leftrightarrow \ g(u)=f(u) \ (a.e.)$. 
\end{tabular}
\end{center}
Let us consider the \gd\ for regression. Suppose that $g(x,y)$, $g(y|x)$, and $g(x)$ are the underlying probability density functions of $(x,y)$, $y$ given $x$, and $x$, respectively. Let $f(y|x)$ be another conditional probability density function of $y$ given $x$. Let $\gamma$ be the positive tuning parameter which controls a trade-off between efficiency and robustness. 
\citet{fujisawa.eguchi.2008} proposed the following cross entropy and divergence:
\begin{align}
& d_{\gamma,1} (g(y|x),f(y|x);g(x))  \nonumber \\
& =  -\frac{1}{\gamma} \log \int  \exp\{ - \gamma d_\gamma(g(y|x),f(y|x)) \}
g(x)dx \nonumber \\
& =  -\frac{1}{\gamma} \log \int \left\{ \int g(y|x) f(y|x)^{\gamma} dy \bigg \slash \left( \int f(y|x)^{1+\gamma}dy\right)^{\frac{\gamma}{1+\gamma}} \right\}g(x)dx \nonumber \\
& =  -\frac{1}{\gamma} \log \int \int \left\{ f(y|x)^{\gamma}  \bigg \slash \left( \int f(y|x)^{1+\gamma}dy\right)^{\frac{\gamma}{1+\gamma}} \right\} g(x,y)dxdy .  \label{type I gam cross} 
\end{align}
\begin{align}
D_{\gamma,1} (g(y|x),f(y|x);g(x)) = -d_{\gamma,1}(g(y|x),g(y|x);g(x))+d_{\gamma,1}(g(y|x),f(y|x);g(x)).   
\end{align}
The cross entropy is empirically estimable, as seen in Section~\ref{Sect: gam-est def}, and the parameter estimation is easily defined. \citet{kawashim.fujisawa.2017} proposed the following cross entropy and divergence:
%
\begin{align}
& d_{\gamma,2} (g(y|x),f(y|x);g(x)) \nonumber \\
 &   = \mathalpha{-} \frac{1}{\gamma} \log \! \int \! \left( \! \int \! g(y|x) f(y|x)^{\gamma} dy  \right) \! g(x) dx \mathalpha{+} \frac{1}{1+\gamma} \log \! \int \! \left( \! \int \! f(y|x)^{1+\gamma} dy  \right) \! g(x) dx \nonumber \\
 &   =  -\frac{1}{\gamma} \log \int  \int  f(y|x)^{\gamma} g(x,y) dxdy + \frac{1}{1+\gamma} \log \int \left( \int f(y|x)^{1+\gamma} dy  \right) g(x) dx.  \label{type II gam cross}
\end{align}
\begin{align}
D_{\gamma,2} (g(y|x),f(y|x);g(x)) = -d_{\gamma, 2}(g(y|x),g(y|x);g(x))+d_{\gamma, 2}(g(y|x),f(y|x);g(x)). 
\end{align}
The base measures on the explanatory variable are taken twice on each term of \gd\ for the \iidp. This extension from the \iidp\ to the regression problem seems to be more natural than \eqref{type I gam cross}. The cross entropy is also empirically estimable. We call these two type I and type II, respectively. These two divergences satisfy the following two basic properties of divergence: 
\begin{center}
\begin{tabular}{cl}
(i) & $D_{\gamma,j}(g(y|x),f(y|x);g(x)) \ge 0$. \\
(ii) & $D_{\gamma,j}(g(y|x),f(y|x);g(x)) = 0 \ \Leftrightarrow \ g(y|x)=f(y|x) \ (a.e.)$. 
\end{tabular}
\end{center}
The equality holds for the conditional \pdf\ instead of usual \pdf. 

%
%

Theoretical properties of $\gamma$-divergence for the \iidp\ were deeply investigated by \citet{fujisawa.eguchi.2008}. Theoretical properties of $\gamma$-divergence for regression were studied by \citet{fujisawa.eguchi.2008}, \citet{kanamori.fujisawa.2015} and \citet{kawashim.fujisawa.2017}, but not well under heterogeneous contamination, which is special in the regression problem and does not appear in the i.i.d. case. \citet{hung.etal.2017} pointed out that a logistic regression model with mislabel can be regarded as a logistic regression model with heterogeneous contamination and then applied the type I to a usual logistic regression model, which enables us to estimate the parameter of the logistic regression model without estimating a mislabel model even if mislabels exist. They also investigated theoretical properties of robustness, but they assumed that $\gamma$ is sufficiently large. In Section~\ref{Sect: rob prop}, we will see that the type I is superior to type II under hetrogeneous contamination in the sense of strong robustness without assuming that $\gamma$ is sufficiently large. Here we mention that the density power divergence \citep{basu.etal.1998} is another candidate of divergence which gives robustness, but it does not have strong robustness \citep{fujisawa.eguchi.2008,hung.etal.2017}.

\subsection{Estimation for $\gamma$-regression}\label{Sect: gam-est def}

Let $f(y|x;\theta)$ be a conditional probability density function of $y$ given $x$ with parameter $\theta$. Let $(x_1,y_1) ,\ldots, (x_n,y_n)$ be the observations randomly drawn from the underlying distribution $g(x,y)$. 
Using the formulae (\ref{type I gam cross}) and (\ref{type II gam cross}), the $\gamma$-cross entropy for regression can be empirically estimated by 
\begin{eqnarray*}
\bar{d}_{\gamma,1} (f(y|x;\theta)) &=&  -\frac{1}{\gamma} \log \frac{1}{n} \sum_{i=1}^n  \frac{   f(y_i|x_i ; \theta)^{\gamma}  }{  \left( \int f(y|x_i ;\theta)^{1+\gamma}dy\right)^\frac{\gamma}{1+\gamma}}, \\ 
\bar{d}_{\gamma,2} (f(y|x;\theta)) &=& - \frac{1}{\gamma} \! \log \! \left\{ \! \frac{1}{n} \! \sum_{i=1}^n  f(y_i | x_i ;\theta)^{\gamma} \! \right\} \mathalpha{+} \frac{1}{1+\gamma} \! \log \! \left\{ \! \frac{1}{n} \! \sum_{i=1}^n \int f(y | x_i ;\theta)^{1+\gamma} dy \! \right\} .
\end{eqnarray*}
The estimator can be defined as the minimizer by
\begin{align*}
\hat{\theta}_{\gamma,j} =\argmin_{\theta} \bar{d}_{\gamma,j}(f(y|x;\theta)) \qquad \mbox{for $j=1,2$}.
\end{align*}
In a similar way to in \citet{fujisawa.eguchi.2008}, we can show that $\hat{\theta}_{\gamma,j}$ converges to $\theta^*_{\gamma,j}$ for $j=1,2$, where
\begin{align*}
\theta^*_{\gamma,j} & = \argmin_{\theta} D_{\gamma,j}(g(y|x),f(y|x;\theta);g(x)) \\
&= \argmin_{\theta} d_{\gamma,j}(g(y|x),f(y|x;\theta);g(x)).
\end{align*}

Suppose that $f(y|x;\theta^*)$ is the target conditional \pdf. The latent bias is expressed as $\theta^*_{\gamma,j}-\theta^*$. This is zero when the underlying model belongs to a parametric model, in other words, $g(y|x)=f(y|x;\theta^*)$, but not always zero when the underlying model is contaminated by outliers. This issue will be discussed in Section~\ref{Sect: rob prop}.

\subsection{Case of location-scale family} \label{Sect: location scale}

Here we show that both types of \gd\ give the same parameter estimation when the parametric conditional probability density function $f(y|x;\theta)$ belongs to a location-scale family in which the scale does not depend on the explanatory variable, given by 
\begin{align}
f(y|x;\theta) = \frac{1}{\sigma} s \left( \frac{y-q(x;\zeta)}{\sigma} \right), \label{eq:location_scale}
\end{align}
where $s(y)$ is a probability density function, $\sigma$ is a scale parameter and $q(x;\zeta)$ is a location function with a regression parameter $\zeta$, e.g., $q(x;\zeta)=x^T \zeta $. 
Then, we can obtain 
\begin{align}\label{location-scale integral}
\int f(y|x;\theta)^{1+\gamma} dy &= \int \frac{1}{\sigma^{1+\gamma}} s \left( \frac{y-q(x;\zeta)}{\sigma} \right)^{1+\gamma} dy \nonumber \\
&= \sigma^{-\gamma} \int s(z)^{1+\gamma} dz.
\end{align}
This does not depend on the explanatory variable $x$.
Using this property, we can show that both types of $\gamma$-cross entropy are the same as follows: 
\begin{align*}
& d_{\gamma,2} (g(y|x),f(y|x;\theta);g(x)) \nonumber \\
& =  -\frac{1}{\gamma} \log \int \left\{ \int g(y|x) f(y|x;\theta)^{\gamma} dy \bigg \slash \left( \int f(y|x;\theta)^{1+\gamma}dy\right)^{\frac{\gamma}{1+\gamma}} \right\}g(x)dx \\
& =  -\frac{1}{\gamma} \log  \left\{ \int \int g(x,y) f(y|x;\theta)^{\gamma} dx dy \bigg \slash \left( \int f(y|x;\theta)^{1+\gamma}dy\right)^{\frac{\gamma}{1+\gamma}} \right\}  \\
& =  -\frac{1}{\gamma} \log \int \int g(x,y) f(y|x;\theta)^{\gamma} dx dy + \frac{1}{1+\gamma}  \log \int f(y|x;\theta)^{1+\gamma}dy  \\
& =  -\frac{1}{\gamma} \log \int \int g(x,y) f(y|x;\theta)^{\gamma} dx dy + \frac{1}{1+\gamma}  \log \int f(y|x;\theta)^{1+\gamma}dy \int g(x) dx  \\
& =  d_{\gamma,2} (g(y|x),f(y|x;\theta);g(x)). 
\end{align*}
The second equality holds from \eqref{location-scale integral}. As a result, both types of \gd\ give the same parameter estimation, because the estimator is defined by the empirical estimation of cross entropy. 
However, it should be noted that both types of $\gamma$-divergence are not the same, because $d_{\gamma,1} (g(y|x),g(y|x);g(x)) \neq d_{\gamma,2} (g(y|x),g(y|x);g(x)) $.

\section{Robust properties}\label{Sect: rob prop}

In this section, we show a distinct difference between two types of $\gamma$-divergence. 

\subsection{Contamination model and basic condition} \label{Sect: contamination and condition}

Let $\delta(y|x)$ be the contamination conditional probability density function related to outliers. Let $\eps(x)$ and $\eps$ denote the outlier ratios which depends on $x$ and does not, respectively. Suppose that the underlying conditional probability density functions under heterogeneous and homogeneous contaminations are given by
\begin{align*}
g(y|x) &= (1-\eps(x))f(y|x;\theta^*) + \eps(x) \delta(y|x), \\
g(y|x) &= (1-\eps)f(y|x;\theta^*) + \eps \delta(y|x).
\end{align*}

Let
\begin{align*}
 {\nu}_{f,\gamma }(x)  = \left\{ \int \delta (y|x) f(y|x)^{ \gamma } dy  \right\}^{ \frac{1}{ \gamma} }, \qquad
{\nu}_{f, \gamma} = \left\{ \int \nu_{f,\gamma}(x)^{\gamma} g(x) dx  \right\}^{ \frac{1}{ \gamma }} .
\end{align*} 
Here we assume that
\begin{align*}
\nu_{f_{\theta^*},\gamma} \approx 0.
\end{align*}
This is an extended assumption used for the \iidp\ \citep{fujisawa.eguchi.2008} to the regression problem. This assumption implies that $\nu_{f_{\theta^*},\gamma}(x) \approx 0$ for any $x$ (a.e.) and  illustrates that the contamination conditional probability density  function $\delta(y|x)$ lies on the tail of the target conditional probability density function $f(y|x;\theta^*)$. 
For example, if $\delta(y|x)$ is the dirac function at the outlier $y_{\dag}(x)$ given $x$, then we have $\nu_{f_{\theta^*},\gamma}(x) = f(y_{\dag}(x)|x;\theta^*) \approx 0$, which is reasonable because $y_{\dag}(x)$ is an outlier. 

Here we also consider the condition $\nu_{f_{\theta},\gamma} \approx 0$, which is used later. This will be true in the neighborhood of $\theta=\theta^*$. In addition, even when $\theta$ is not close to $\theta^*$, if $\delta(y|x)$ lies on the tail of $f(y|x;\theta)$, we cansee $\nu_{f_{\theta},\gamma} \approx 0$. 

To make the discussion easier, we prepare the monotone transformation of both types of $\gamma$-cross entropies for regression by
\begin{align*}
\tilde{d}_{\gamma,1}(g(y|x),f(y|x;\theta);g(x)) &= - \exp  \left\{ -\gamma d_{\gamma,1}(g(y|x),f(y|x;\theta);g(x)) \right\} \\
&= - \int \int \frac{   f(y|x;\theta)^{\gamma}  }{ \left(  \int f(y|x;\theta)^{1+\gamma} dy  \right)^{\frac{\gamma}{1+\gamma} } } g(y|x) g(x)dx dy  , \\
\tilde{d}_{\gamma,2}(g(y|x),f(y|x;\theta);g(x)) &= - \exp  \left\{ -\gamma d_{\gamma,2}(g(y|x),f(y|x;\theta);g(x)) \right\} \\
&=  - \frac{  \int \left(  \int g(y|x) f(y|x;\theta)^{\gamma} dy \right) g(x) dx    }{ \left\{ \int \left( \int f(y|x;\theta)^{1+\gamma} dy  \right) g(x) dx \right\}^{\frac{\gamma}{1+\gamma} } } . 
\end{align*}

\subsection{Type-I of $\gamma$-divergence}

We see
\begin{eqnarray*}
\lefteqn{ \tilde{d}_{\gamma,1}(g(y|x),f(y|x;\theta);g(x)) } \\
 &=& - \int \frac{ \int g(y|x) f(y|x;\theta)^{\gamma} dy  }{ \left( \int f(y|x;\theta)^{1+\gamma} dy  \right)^{\frac{\gamma}{1+\gamma}} } g(x) dx \\
 &=& - \int \frac{  \int \left\{ (1-\eps(x)) f(y|x;\theta^*) + \eps(x) \delta(y|x) \right\} f(y|x;\theta)^{\gamma} dy }{ \left( \int f(y|x;\theta)^{1+\gamma} dy  \right)^{\frac{\gamma}{1+\gamma}} } g(x) dx \\
 &=& - \int \frac{ \int  f(y|x;\theta^*) f(y|x;\theta)^{\gamma} dy }{ \left( \int f(y|x;\theta)^{1+\gamma} dy  \right)^{\frac{\gamma}{1+\gamma} }} (1-\eps(x)) g(x) dx - \int \frac{  \int  \delta(y|x;\theta) f(y|x;\theta)^{\gamma} dy  }{ \left( \int f(y|x;\theta)^{1+\gamma} dy  \right)^{\frac{\gamma}{1+\gamma}} } \eps(x) g(x) dx \\
 &=& - \tilde{d}_{\gamma,1}(f(y|x;\theta^*),f(y|x;\theta)); \tilde{g}(x)) - \int \frac{ \nu_{f_{\theta},\gamma}(x)^{\gamma}  }{ \left(   \int f(y|x;\theta)^{1+\gamma} dy \right)^{\frac{\gamma}{1+\gamma} } }   \eps(x) g(x) dx   , \label{heter prop eq type2}
\end{eqnarray*}
where $\tilde{g}(x) = (1-\eps(x))g(x)$. From this relation, we can easily show the following theorem. 

\begin{theorem}\label{hetero_cross type1}
Under the condition $\nu_{f_{\theta},\gamma} \approx 0$ and $\int f(y|x;\theta)^{1+\gamma} dy > 0 $, we have
\begin{align*}
\tilde{d}_{\gamma,1}(g(y|x),f(y|x;\theta);g(x))  \approx  \tilde{d}_{\gamma,1} (f(y|x;\theta^*),f(y|x;\theta);  \tilde{g}(x)) .
\end{align*}
\end{theorem}

Using this theorem, we can expect that the latent bias $\theta^*_{\gamma,1}-\theta^*$ is close to zero, because 
\begin{eqnarray*}
\arg\min_\theta \tilde{d}_{\gamma,1}(g(y|x),f(y|x;\theta);g(x)) &=& \arg\min_\theta d_{\gamma,1}(g(y|x),f(y|x;\theta);g(x)) \,=\, 
\theta^*_{\gamma,1} \\
\arg\min_\theta \tilde{d}_{\gamma,1} (f(y|x;\theta^*),f(y|x;\theta);  \tilde{g}(x)) &=& 
\arg\min_\theta d_{\gamma,1} (f(y|x;\theta^*),f(y|x;\theta);  \tilde{g}(x)) \,=\,
 \theta^*. 
\end{eqnarray*}
The last equality holds even when $g(x)$ is replace by $\tilde{g}(x) = (1-\eps(x))g(x)$. 

In addition, we can have the modified Pythagorian relation approximately. 

\begin{theorem}\label{hetero_pytha type1}
Under the condition $\nu_{f_{\theta},\gamma} \approx 0$ and $\int f(y|x;\theta)^{1+\gamma} dy > 0 $, the modified Pythagorian relation among $g(y|x)$, $f(y|x;\theta^*)$, $f(y|x;\theta)$ approximately holds:
$$ D_{\gamma,1}(g(y|x),f(y|x;\theta);g(x)) \approx D_{\gamma,1}(g(y|x),f(y|x;\theta^*);g(x)) + D_{\gamma,1}(f(y|x;\theta^*),f(y|x;\theta);\tilde{g}(x)). $$
\end{theorem}

The Pythagorian relation implies that the minimizer of $D_{\gamma,1}(g(y|x),f(y|x;\theta);g(x))$ is almost the same as the minimizer of $D_{\gamma,1}(f(y|x;\theta^*),f(y|x;\theta);\tilde{g}(x))$, which is $\theta^*$. This also implies the strong robustness. 

In the theorems, we assume $\nu_{f_{\theta},\gamma} \approx 0$ and $\int f(y|x;\theta)^{1+\gamma} dy > 0 $. The former condition was already discussed in Section~\ref{Sect: contamination and condition}. Here we investigate the latter condition. When the parametric conditional \pdf\ belongs to a location-scale family \eqref{eq:location_scale}, this condition will be expected to hold, because
\begin{eqnarray*}
\int f(y|x;\theta)^{1+\gamma} dy
 &=& \int \frac{1}{\sigma^{1+\gamma}} s \left( \frac{y-q(x;\zeta)}{\sigma} \right)^{1+\gamma} dy = \frac{1}{\sigma^{\gamma}} \int s \left( z \right)^{1+\gamma} dz. 
\end{eqnarray*}
We can also verify that this condition holds for a logistic regression model, a Poisson regression model, and so on. 

Finally we mention the homogeneous contamination. The modified Pythagorian relation in Theorem~\ref{hetero_pytha type1} is changed to the usual Pythagorian relation, because we can easily see $D_{\gamma,1}(f(y|x;\theta^*),f(y|x;\theta);\tilde{g}(x))=D_{\gamma,1}(f(y|x;\theta^*),f(y|x;\theta);g(x))$ under homogeneous contamination.

\subsection{Type 2 of $\gamma$-divergence}

First, we illustrate that the strong robustness does not hold in general under heterogeneous contamination, unlike for type 1. We see
\begin{eqnarray*}
\lefteqn{ \tilde{d}_{\gamma,2}(g(y|x),f(y|x;\theta);g(x)) } \\
 &=& - \frac{  \int \left(  \int g(y|x) f(y|x;\theta)^{\gamma} dy \right) g(x) dx    }{ \left\{ \int \left( \int f(y|x;\theta)^{1+\gamma} dy  \right) g(x) dx \right\}^{\frac{\gamma}{1+\gamma} } } \\
 &=& - \frac{  \int \left(  \int  (1-\eps(x)) f(y|x;\theta^*) f(y|x;\theta)^{\gamma} dy + \int  \eps(x) \delta(y|x) f(y|x;\theta)^{\gamma} dy \right) g(x) dx }{ \left\{ \int \left( \int f(y|x;\theta)^{1+\gamma} dy  \right) g(x) dx \right\}^{\frac{\gamma}{1+\gamma} } } \\
 &=& - \frac{  \int \left(  \int  (1-\eps(x)) f(y|x;\theta^*) f(y|x;\theta)^{\gamma} dy + \int  \eps(x) \nu_{f_\theta,\gamma}(x) \right) g(x) dx }{ \left\{ \int \left( \int f(y|x;\theta)^{1+\gamma} dy  \right) g(x) dx \right\}^{\frac{\gamma}{1+\gamma} } } \\
 & \approx & - \frac{  \int \int  f(y|x;\theta^*) f(y|x;\theta)^{\gamma} dy (1-\eps(x)) g(x) dx }{ \left\{ \int \left( \int f(y|x;\theta)^{1+\gamma} dy  \right) g(x) dx \right\}^{\frac{\gamma}{1+\gamma} } }.
\end{eqnarray*}
The last aproximation holds from $\nu_{f_\theta,\gamma}(x) \approx 0$. This can not be expressed using $d_\gamma(f(y|x;\theta^*),f(y|x;\theta);h(x))$ with an appropriate base measure $h(x)$, unlike for type 1, because the base measure of the numerator on the explanatory variable is different from that of the denominator. As in numerical experiments, the type 2 presents a significant bias under heterogenous contamination. However, as already mentioned, when the parametric conditional \pdf\ belongs to a location-scale family \eqref{eq:location_scale}, the cross entropy for type 2 is the same as that for type 1 and then the type 2 has the strong robustness. In addition, under homogeneous contamination, we have $\tilde{d}_{\gamma,2}(g(y|x),f(y|x;\theta);g(x)) \approx (1-\eps) \tilde{d}_{\gamma,2}(f(y|x;\theta^*),f(y|x;\theta);g(x))$ and then we expect that the latent bias $\theta^*_{\gamma,2}-\theta^*$ is sufficiently small.

\section{Numerical experiment}\label{Sect: experment}

In this section, using a simulation model, we compare the type 1 with the type 2. 


As shown in Section~\ref{Sect: rob prop}, the distinct difference occurs under heterogeneous contamination when the parametric conditional probability density function $f(y|x;\theta)$ does not belong to a location-scale family. 
Therefore, we used the logistic regression model as the simulation model, given by 
\begin{align*}
Pr(y=1|x) = \pi(x;\beta) , \ Pr(y=0|x) = 1-\pi(x;\beta), 
\end{align*}
where $\pi(x;\beta)= \{1+\exp (- \beta_{0} - x_1 \beta_1- \cdots - x_p \beta_p  )\}^{-1}$.
%
%
%
%
The sample size and the number of explanatory variables were set to be $n=1000$ and $p=5$, respectively. 
The true coefficients were given by 
\begin{gather*}
\beta_{0}=0, \ \beta_{1} =1 ,\ \beta_{2} = -1,\ \beta_{3} = 1,\ \beta_{4} =-1 ,\ \beta_{5} =0.  
\end{gather*}
%
%
The explanatory variables were generated from a normal distribution $N(0,\Sigma)$ with $ \Sigma=(0.2^{ | i-j | })_{1 \leq i,j \leq p }$. 
We generated 100 random samples. 

Outliers were incorporated into simulations. 
We investigated four outlier ratios ($\eps=0.1, \ 0.2, \ 0.3 \mbox{ and } 0.4$) and the following outlier pattern: The outliers were generated around the edge part of the explanatory variable, where the explanatory variables were generated from $N(\bm{\mu}_{\rm out}, 0.5^{2} {\bf I})$ where ${\bm{\mu}_{\rm out}}=(20,0,20,0,0)$ and the response variable $y$ is set to $0$. 
%

In order to verify the fitness of regression coefficient, we used the mean squared error (MSE) as the performance measure, given by 
\begin{align*}
\textrm{MSE} &= \frac{1}{p+1} \sum_{j=0}^p (\hat{\beta}_j - \beta_j^* )^2, 
\end{align*}
where $\beta_j^*$'s are the true coefficients. 
The tuning parameter $\gamma$ in the $\gamma$-divergence was set to $0.5$ and $1.0$. 

Table 1 shows the MSE in the case $\eps=0.1$, $0.2$, $0.3$ and $0.4$.
The type 2 presented smaller MSEs than the type 1. The difference between two types was larger as the outlier ratio $\eps$ was larger. 

\begin{table}[H]
 \centering
\caption{MSE under heterogeneous contamination}
\vspace*{3mm}
\begin{tabular}{|c|cc|} \hline
Methods  &  $\gamma=0.5$ & $\gamma=1.0$  \\ \hline
 & \multicolumn{2}{|c|}{$\eps=0.1$} \\ \hline
Type 1 & 0.00620 & 0.00712  \\
Type 2 & 0.00810 & 0.0276  \\ \hline
 & \multicolumn{2}{|c|}{$\eps=0.2$} \\ \hline
Type 1 & 0.0136 & 0.0149  \\
Type 2 & 0.0215 & 0.110  \\ \hline
 & \multicolumn{2}{|c|}{$\eps=0.3$} \\  \hline
Type 1 & 0.0262 & 0.0282  \\
Type 2 & 0.0472 & 0.282  \\ \hline
 & \multicolumn{2}{|c|}{$\eps=0.4$} \\ \hline
Type 1 & 0.0514 & 0.0547  \\
Type 2 & 0.0998 & 0.648  \\ \hline
\end{tabular}
\end{table}

\bibliography{ref}

\end{document}